# A global mathematical investigation of a predator-prey model

S. A. Treskov*, E. P. Volokitin†

**Abstract**

We construct a global bifurcation diagram of the plane differential system

$$\dot{x} = x(1-x) - xy/(a+x^2),\ \dot{y} = y(\delta - \beta y/x),$$
$$x(t) > 0, y(t) > 0, a > 0, \delta > 0, \beta > 0,$$

which describes the predator-prey interaction.

In [1] was considered the following planar differential systems which describes a mathematical predator-prey model

$$\begin{aligned} \dot{x} &= rx(1-\tfrac{x}{K}) - \tfrac{mx}{b+x^2}y, \\ \dot{y} &= sy(1-\tfrac{y}{hx}). \end{aligned} \tag{1}$$

Here $x(t) > 0,\ y(t) > 0$ — the population densities of the prey and the predator respectively and $r, s, K, h, m$ — positive parameters ( for more details see [1]).

After the variables and parameters scaling the system takes the form

$$\begin{aligned} &\dot{x} = x(1-x) - \tfrac{xy}{a+x^2}, \\ &\dot{y} = y(\delta - \beta\tfrac{y}{x}), \\ &x(t) > 0,\ y(t) > 0,\ a > 0,\ \delta > 0,\ \beta > 0. \end{aligned} \tag{2}$$

System (2) was investigated with the help of the qualitative theory of the ordinary differential equations in [1] where it was demonstrated some properties of the system. For example, it was shown that there is an positively invariant region in the phase plane and it was found how many steady states are located there and it were found their types and so on. During the investigation the authors of [1] studied some bifurcations of codim 1 and 2: the turning point bifurcation, the cusp bifurcation, the Andronov-Hopf bifurcation, the Bogdanov-Takens bifurcation, the loop bifurcation. A special attention was given to the study of periodic solutions. The result of [1] were extended in [2].

*The work is supported by RFBR (09-01-00070) and by Siberian branch of RAS (interdisciplinaire integration projects no. 107, 119

†The work is supported by RFBR (09-01-00070) and by Siberian branch of RAS (interdisciplinaire integration projects no. 107, 119

In [3] we derived some results which supplemented the investigation from [1, 2].

In this paper we fulfil a complete parametric analysis of system (2) taking into account mentioned above and our new results. In particular, we consider all bifurcations of system (2) of codim 1, 2 and 3. The result of our study is the list of all rough phase portraits of the system under study and we can divide the entire parameter space $R^3_+$ into regions with different rough phase portraits. Our results give the base for detailed simulations.

**1.** For any parameter values $a, \beta, \delta$ the domain $\Pi = \{(x, y) \in R^2 : 0 < x < 1, y > 0\}$ is $\omega$-invariant set of (2). The boundary of the domain $\Pi$ contains two steady states of (2): points $(0, 0)$ and $(1, 0)$.

For all parameter values the rest point $(1, 0)$ is a hyperbolic saddle with its outgoing separatrix passing into $\Pi$.

The origin $(0, 0)$ is a non-analytic steady state. we change the time $t = \tau x$. After that (2) takes the equivalent form

$$\begin{aligned} &\dot{x} = x^2(1-x) - \tfrac{x^2 y}{a+x^2} \equiv P_2(x,y) + \ldots, \\ &\dot{y} = \delta xy - \beta y^2 \equiv Q_2(x,y), \\ &x(t) > 0, y(t) > 0, a > 0, \delta > 0, \beta > 0. \end{aligned} \tag{3}$$

We denote the derivation on $\tau$ with the point again. $P_2(x, y), Q_2(x, y)$ are homogeneous polynomials of degree two, dots replace terms of degree more than two.

For (3) the origin is a complex analytic steady state, the matrix of the linear part contains zero elements only.

As it is proved in [4] in our case we have the following

1) if $\delta \leq 1$ the neighborhoods of the origin inside $R^2_+$ contains a unique hyperbolic sector with the coordinate axes as its boundary;

2) if $\delta > 1$ the neighborhoods of the origin inside $R^2_+$ contains one hyperbolic sector and one parabolic sector with a common boundary which is the trajectory outgoing from the origin with the angle $\vartheta = atan\frac{\delta-1}{\beta}$. The second boundary of the parabolic sector is the x-axis. 0-trajectories have the asymptotic $y = kx^\delta, k \geq 0$.

Phase portraits of system (2) near the origin and the infinity are shown in the Fig. 1.

**2.** To fulfil the complete parametric analysis we first need to study all codimension one bifurcations in system (2) to which correspond hypersurfaces in the parametric space dividing the parametri space into regions within which the system has distinct rough phase portraits. For parameter values on this surfaces the system has one from the following singularities

I. The turning point (the saddle-node).

II. The fine focus of the first degree.

III. The double limit cycle.

IV. The seaparatrix loop with the non-zero saddle value.

The local bifurcations of codim 1 for I, II can be described by explicit formulas

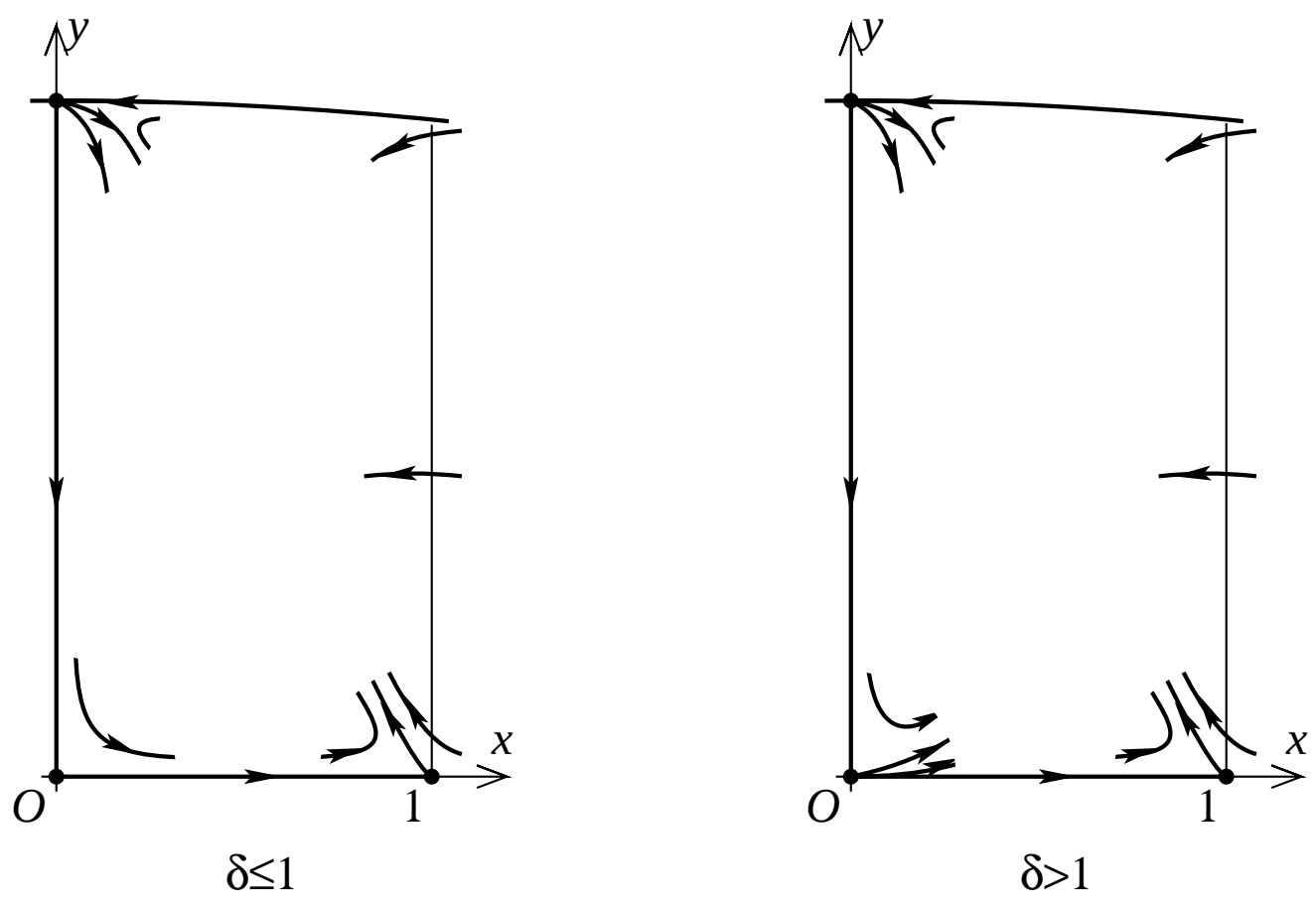

Fig. 1.

Let $(x_0, y_0)$ is a steady state of (2). Then the following algebraic system is valid

$$x_0^3 - x_0^2 + (a + \delta/\beta)x_0 - a = 0, \quad y_0 = (\delta/\beta)x_0,$$

which has from one to three solutions inside our domain $\Pi$.

Let $J$ is a matrix of a linear part near the steady state. We have

$$F \equiv x_0^3 - x_0^2 + (a + \delta/\beta)x_0 - a = 0, \; y_0 = \tfrac{\delta}{\beta}x_0,$$
$$\det J = k_1(-a^2\beta + 2a^2\beta x_0 + 2a\delta x_0 - 2a\beta x_0^2 + 4a\beta x_0^3 - \beta x_0^4 + 2\beta x_0^5),$$
$$\mathrm{tr} J = k_2(a\delta + ax_0 - 2x_0^2 + \delta x_0^2 + 3x_0^3), \; k_1 > 0, \; k_2 > 0.$$

I. After the elimination of $x$ conditions $F = 0, \det J = 0$ isolate the bifurcation set $T$ in the parameter space. The set corresponds to the double steady state. $T$ is a piecewise smooth surface and it may be defined as follows

$$a = \tau^2 - 2\tau^3, \; \beta = \frac{\delta}{2\tau(1-\tau)^2}.$$

The surface is sketched in the Fig. 2. The points inside the angle correspond to the existence of three nonhyperbolic steady states in system (2) (one saddle and two antisaddles) and the points outside the angle correspond to the existence of a unique steady state (the node or the focus). The points at the surface $T$ (but not at $TT$) correspond to the existence of one double saddle-node and one node or focus. The points at the curve $TT$ correspond to the existence of one steady state of the multiplicity 3.

Here and below we examinated the corresponding non-degenerate conditions of the bifurcations with the help of results from [6], [2] and our results.

II. After the elimination of $x_0$ conditions $F = 0, \mathrm{tr} J = 0, \det J > 0$ isolate the bifurcation set $H$ corresponding to the existence of the fine focus in system

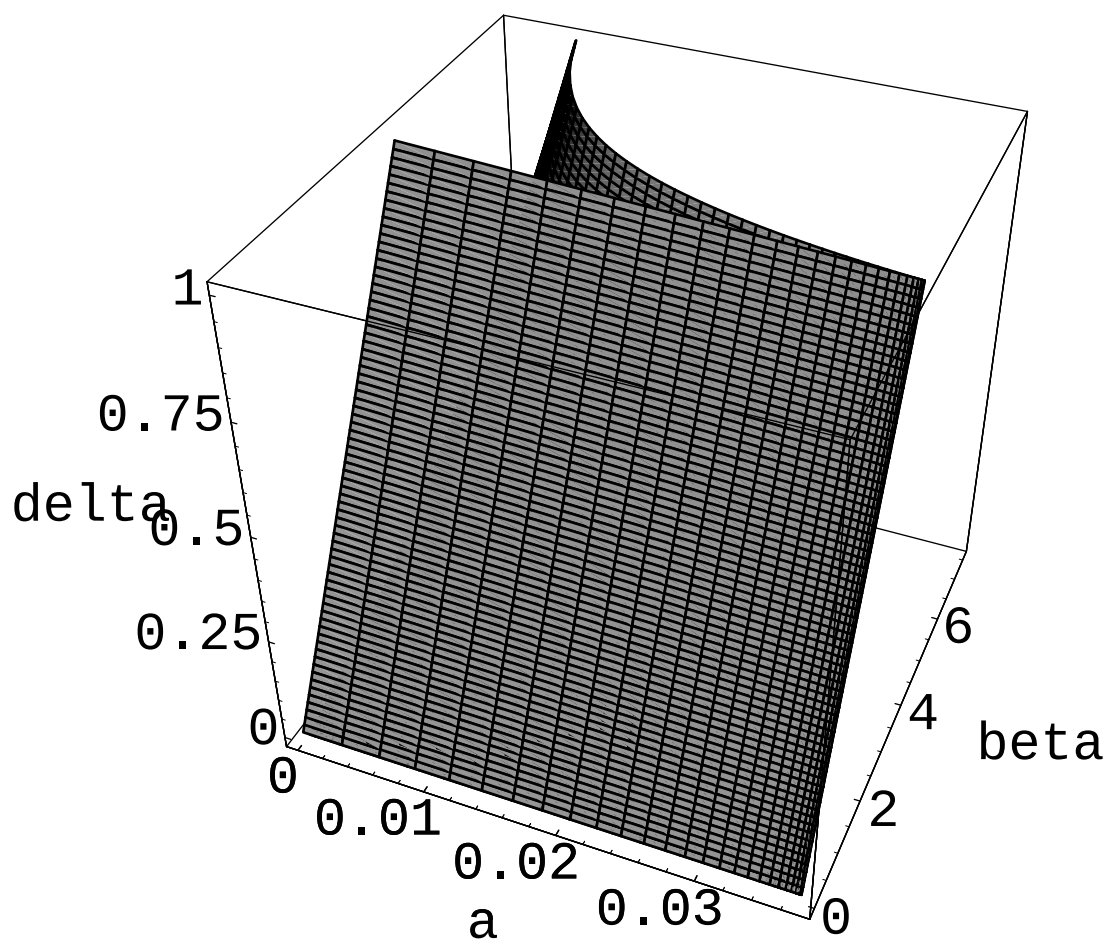

Fig. 2.

(2) (the eigenvalues of the matrix of a linear part are pure imaginary). This bifurcation set is a smooth surface with the boundary ($\det J = 0$). It may be defined as follows

$$a = \frac{\tau^2(2-\delta-3\tau)}{\delta+\tau},\ \beta = \frac{\delta(\delta+\tau)}{2\tau(1-\tau)^2}, 1-\delta-\tau > 0.$$

At all cases the steady state $(x_0, y_0)$ is given by the expressions

$$x_0 = \tau,\ y_0 = \frac{\delta}{\beta}\tau.$$

Surface $H$ has a complicate structure. Therefore we do not sketch it here.[1]

In Fig. 3 it is shown schematically sections of surfaces $H$, $T$ with the hyperplanes $\delta = c$ for various values of $c$. On this sections the surfaces are shown as the curves which we label as $H, T$ again. About other notations at the picture we explain below.

Curve $T$ has a cusp. Dy dotted line it is shown the section of an extension of surface $H$ for which $F = 0, \mathrm{tr}J = 0, \det J < 0$. For corresponding parameter values there is a saddle with zero sum of eigenvalues of matrix $J$.

In the Fig. 4 it is shown the section of surfaces $H, T$ by the hyperplane $\delta = .25$. We see that curve $H$ has the selfintersection point which is in the domain where system (2) has three steady states.

The conditions about double cycles (III) and separatrix loops (IV) can not be describe by explicit formulas in our case. To construct surfaces of this nonlocal bifurcations we consider local bifurcations of more high codimensions (two and three). Such a local bifurcations are listed in the table.

[1]In [7] there is a picture with Andronov-Hopf bifurcation surface which is analogous to our surface $H$.

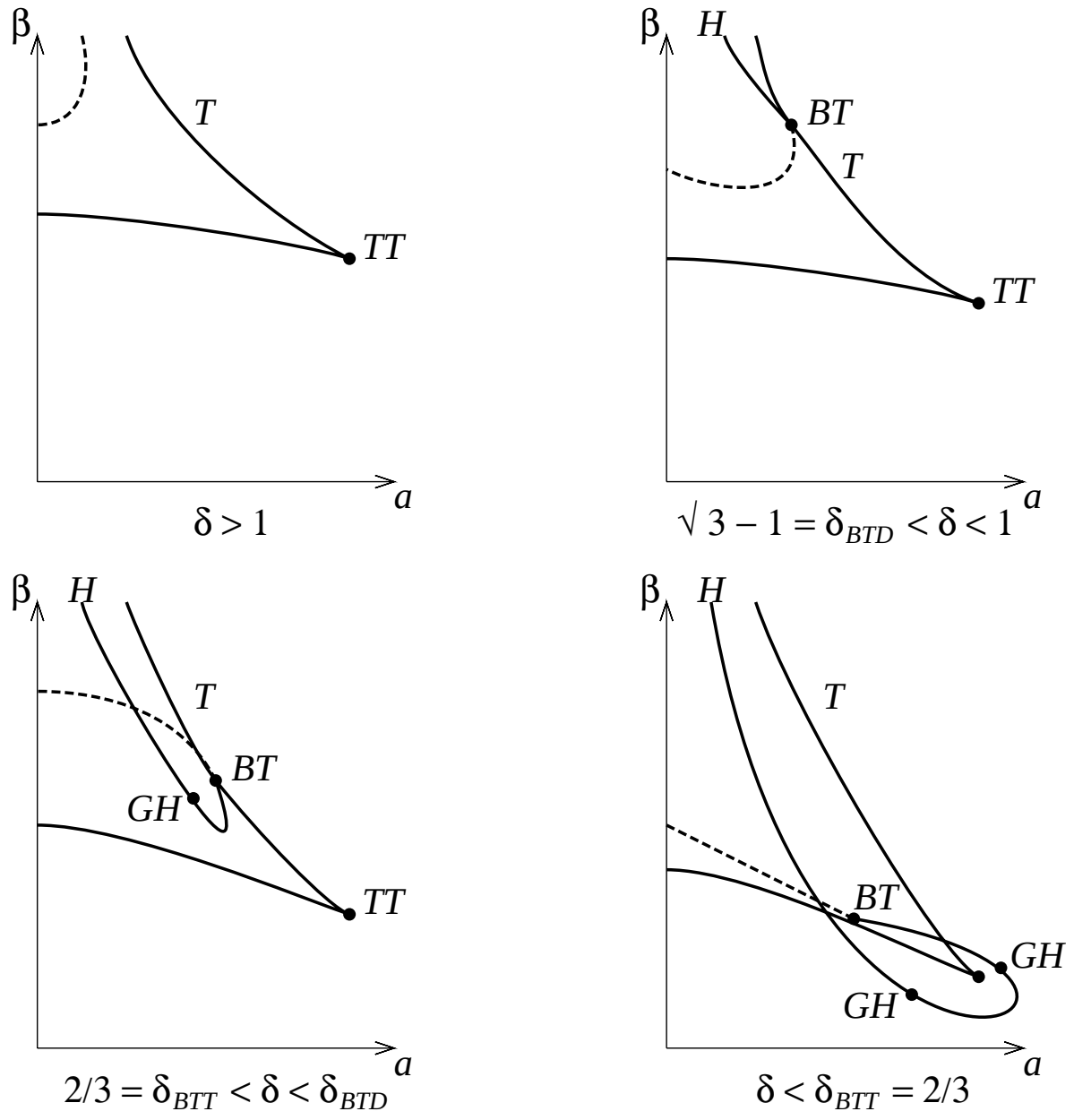

Fig. 3.

Table 1: The list of local bifurcations

| Bifurcation | Notation |
|---|---|
| codim 1 | |
| The turning point | $T$ |
| Andronov-Hopf bifurcation | $H$ |
| codim 2 | |
| The cusp | $TT$ |
| Bogdanov-Takens bifurcation | $HT$ |
| The generalized Hopf bifurcation | $GH$ |
| codim 3 | |
| The focus singularity with nilpotent linear part | $HTT$ |
| The cusp singularity with nilpotent linear part | $DHT$ |

The listed local bifurcations may be described by explicit analitic expressions.

If there is a steady state of multiplicity 3 in (2) then we have a solution of equation $F = 0$ with the same multiplicity. In the parameter space corresponding set is curve $TT$ at the surface $T$ which is given as follows

$$a = 1/27,\ \beta = 27\delta/8,\ \delta > 0.$$

The bifurcation set $HT$ is given by the condition $\lambda_1 = \lambda_2 = 0$ ($F = \det J = \mathrm{tr} J = 0$). The set is a smooth curve and may be described as follows

$$a = \tau^2 - 2\tau^3,\ \beta = \frac{1}{2\tau(1-\tau)},\ \delta = 1 - \tau.$$

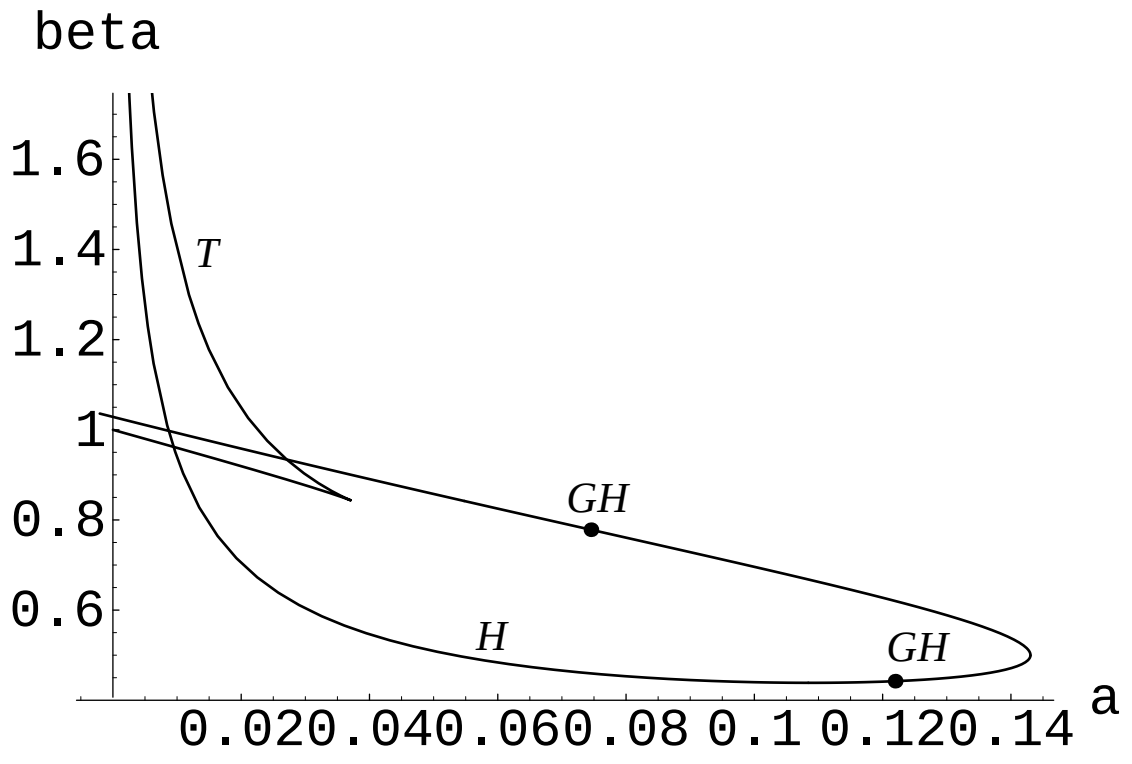

Fig. 4.

The bifurcation set $GH$ is given by the condition $F = \mathrm{tr}J = l_1 = 0,\ \det J > 0$, where $l_1$ is the first Lyapunov value of the fine focus. The non-degenerate condition of this codim 2 bifurcation is the condition $l_2 \neq 0$ where $l_2$ is the second Lyapunov value.

With the help of the algorithm of Poincaré [9, 10] we derived symbolic expressions for the first and second Lyapunov values of the system

$$\dot{x} = f(x, y),\ \dot{y} = g(x, y)$$

These expressions are

$$l_1 = -a_{02}a_{11} + a_{12} - a_{11}a_{20} + 3a_{30} - 2a_{02}b_{02} + 3b_{03} + b_{02}b_{11} + 2a_{20}b_{20} + b_{11}b_{20} + b_{21},$$

$l_2 = 10a_{02}{}^3a_{11} + 30a_{02}a_{03}a_{11} - 9a_{04}a_{11} - 10a_{02}{}^2a_{12} - 15a_{03}a_{12} - 15a_{02}a_{13}+$
$9a_{14} + 97a_{02}{}^2a_{11}a_{20} + 27a_{03}a_{11}a_{20} - 72a_{02}a_{12}a_{20} - 21a_{13}a_{20} + 176a_{02}a_{11}a_{20}{}^2$
$-80a_{12}a_{20}{}^2 + 89a_{11}a_{20}{}^3 + 15a_{02}a_{11}a_{21} - 12a_{12}a_{21} + 18a_{11}a_{20}a_{21} - 3a_{11}a_{22}-$
$75a_{02}{}^2a_{30} - 18a_{03}a_{30} - 6a_{11}{}^2a_{30} - 288a_{02}a_{20}a_{30} - 267a_{20}{}^2a_{30} - 27a_{21}a_{30}-$
$9a_{02}a_{31} - 27a_{20}a_{31} + 9a_{32} + 15a_{11}a_{40} + 45a_{50} + 20a_{02}{}^3b_{02} + 60a_{02}a_{03}b_{02}-$
$18a_{04}b_{02} + 34a_{02}a_{11}{}^2b_{02} - 28a_{11}a_{12}b_{02} + 174a_{02}{}^2a_{20}b_{02} + 24a_{03}a_{20}b_{02}+$
$31a_{11}{}^2a_{20}b_{02} + 208a_{02}a_{20}{}^2b_{02} + 18a_{20}{}^3b_{02} - 30a_{20}a_{21}b_{02} + 12a_{22}b_{02}-$
$63a_{11}a_{30}b_{02} + 66a_{40}b_{02} + 121a_{02}a_{11}b_{02}{}^2 - 71a_{12}b_{02}{}^2 + 68a_{11}a_{20}b_{02}{}^2-$
$174a_{30}b_{02}{}^2 + 106a_{02}b_{02}{}^3 - 18a_{20}b_{02}{}^3 - 30a_{02}{}^2b_{03} - 45a_{03}b_{03} - 9a_{11}{}^2b_{03}-$
$246a_{02}a_{20}b_{03} - 282a_{20}{}^2b_{03} - 18a_{21}b_{03} - 75a_{11}b_{02}b_{03} - 159b_{02}{}^2b_{03} - 60a_{02}b_{04}-$
$66a_{20}b_{04} + 45b_{05} + 31a_{02}{}^2a_{11}b_{11} + 6a_{03}a_{11}b_{11} - 31a_{02}a_{12}b_{11} - 6a_{13}b_{11}+$
$110a_{02}a_{11}a_{20}b_{11} - 70a_{12}a_{20}b_{11} + 79a_{11}a_{20}{}^2b_{11} - 120a_{02}a_{30}b_{11} - 237a_{20}a_{30}b_{11}+$
$52a_{02}{}^2b_{02}b_{11} - 3a_{03}b_{02}b_{11} - 3a_{11}{}^2b_{02}b_{11} + 56a_{02}a_{20}b_{02}b_{11} - 104a_{20}{}^2b_{02}b_{11}-$
$6a_{21}b_{02}b_{11} - 25a_{11}b_{02}{}^2b_{11} - 53b_{02}{}^3b_{11} - 108a_{02}b_{03}b_{11} - 225a_{20}b_{03}b_{11}-$
$15b_{04}b_{11} + 21a_{02}a_{11}b_{11}{}^2 - 21a_{12}b_{11}{}^2 + 21a_{11}a_{20}b_{11}{}^2 - 63a_{30}b_{11}{}^2 - 4a_{02}b_{02}b_{11}{}^2-$
$85a_{20}b_{02}b_{11}{}^2 - 60b_{03}b_{11}{}^2 - 18b_{02}b_{11}{}^3 + 15a_{02}a_{11}b_{12} - 15a_{12}b_{12} + 12a_{11}a_{20}b_{12}-$
$36a_{30}b_{12} - 30a_{20}b_{02}b_{12} - 27b_{03}b_{12} + 27b_{02}b_{13} + 25a_{02}a_{11}{}^2b_{20} - 25a_{11}a_{12}b_{20}-$
$50a_{02}{}^2a_{20}b_{20} - 12a_{03}a_{20}b_{20} + 22a_{11}{}^2a_{20}b_{20} - 192a_{02}a_{20}{}^2b_{20} - 178a_{20}{}^3b_{20}-$
$6a_{02}a_{21}b_{20} - 36a_{20}a_{21}b_{20} + 6a_{22}b_{20} - 54a_{11}a_{30}b_{20} + 60a_{40}b_{20} + 110a_{02}a_{11}b_{02}b_{20}-$
$66a_{12}b_{02}b_{20} + 16a_{11}a_{20}b_{02}b_{20} - 138a_{30}b_{02}b_{20} + 120a_{02}b_{02}{}^2b_{20} - 136a_{20}b_{02}{}^2b_{20}-$
$66a_{11}b_{03}b_{20} - 180b_{02}b_{03}b_{20} - 25a_{02}{}^2b_{11}b_{20} - 6a_{03}b_{11}b_{20} - 3a_{11}{}^2b_{11}b_{20}-$
$182a_{02}a_{20}b_{11}b_{20} - 265a_{20}{}^2b_{11}b_{20} - 3a_{21}b_{11}b_{20} - 38a_{11}b_{02}b_{11}b_{20} - 104b_{02}{}^2b_{11}b_{20}-$
$43a_{02}b_{11}{}^2b_{20} - 124a_{20}b_{11}{}^2b_{20} - 18b_{11}{}^3b_{20} - 6a_{02}b_{12}b_{20} - 36a_{20}b_{12}b_{20} - 3b_{11}b_{12}b_{20}+$
$9b_{13}b_{20} + 25a_{02}a_{11}b_{20}{}^2 - 25a_{12}b_{20}{}^2 - 16a_{11}a_{20}b_{20}{}^2 - 30a_{30}b_{20}{}^2 + 50a_{02}b_{02}b_{20}{}^2-$
$102a_{20}b_{02}b_{20}{}^2 - 75b_{03}b_{20}{}^2 - 13a_{11}b_{11}b_{20}{}^2 - 61b_{02}b_{11}b_{20}{}^2 - 20a_{20}b_{20}{}^3 - 10b_{11}b_{20}{}^3-$
$25a_{02}{}^2b_{21} - 6a_{03}b_{21} - 3a_{11}{}^2b_{21} - 102a_{02}a_{20}b_{21} - 107a_{20}{}^2b_{21} - 3a_{21}b_{21}-$
$16a_{11}b_{02}b_{21} - 44b_{02}{}^2b_{21} - 43a_{02}b_{11}b_{21} - 82a_{20}b_{11}b_{21} - 18b_{11}{}^2b_{21} - 6b_{12}b_{21}-$
$13a_{11}b_{20}b_{21} - 36b_{02}b_{20}b_{21} - 10b_{20}{}^2b_{21} - 6a_{02}b_{22} - 12a_{20}b_{22} + 3b_{11}b_{22} + 9b_{23}+$
$12a_{02}a_{11}b_{30} - 12a_{12}b_{30} + 15a_{11}a_{20}b_{30} - 9a_{30}b_{30} + 24a_{02}b_{02}b_{30} + 24a_{20}b_{02}b_{30}-$
$36b_{03}b_{30} + 6a_{11}b_{11}b_{30} + 9b_{02}b_{11}b_{30} + 24a_{20}b_{20}b_{30} + 12b_{11}b_{20}b_{30} - 3b_{21}b_{30}+$
$6a_{11}b_{31} + 21b_{02}b_{31} + 15b_{20}b_{31} + 18a_{20}b_{40} + 9b_{11}b_{40} + 9b_{41}.$

Let us note that our expressions coinside (with an accuracy to a positive factor) with expressions derived by other authors (see [5] for example).

In our case we have for system (2) (with an accuracy to a positive factor again)

$$\begin{aligned} l_1 = \; & 6\delta-11\delta^2+8\delta^3-2\delta^4+4x_0-39\delta x_0+42\delta^2 x_0-14\delta^3 x_0-24x_0^2+ \\ & +72\delta x_0^2-33\delta^2 x_0^2++36x_0^3-33\delta x_0^3-12x_0^4. \end{aligned}$$

In Fig. 5 it is shown the curve $l_1 = 0$ and the line $x_0 + \delta = 1$ (the last one describes the boundary of surface H). We see that bifurcation set $GH$ is a smooth two-component curve. The second Lyapunov value is negative in all points from $GH$ (we tested this fact by numerical calculations). Therefor the complex focus in system (2) has its multiplicity no more than two and the focus is stable if its multiplicity is equal to two.

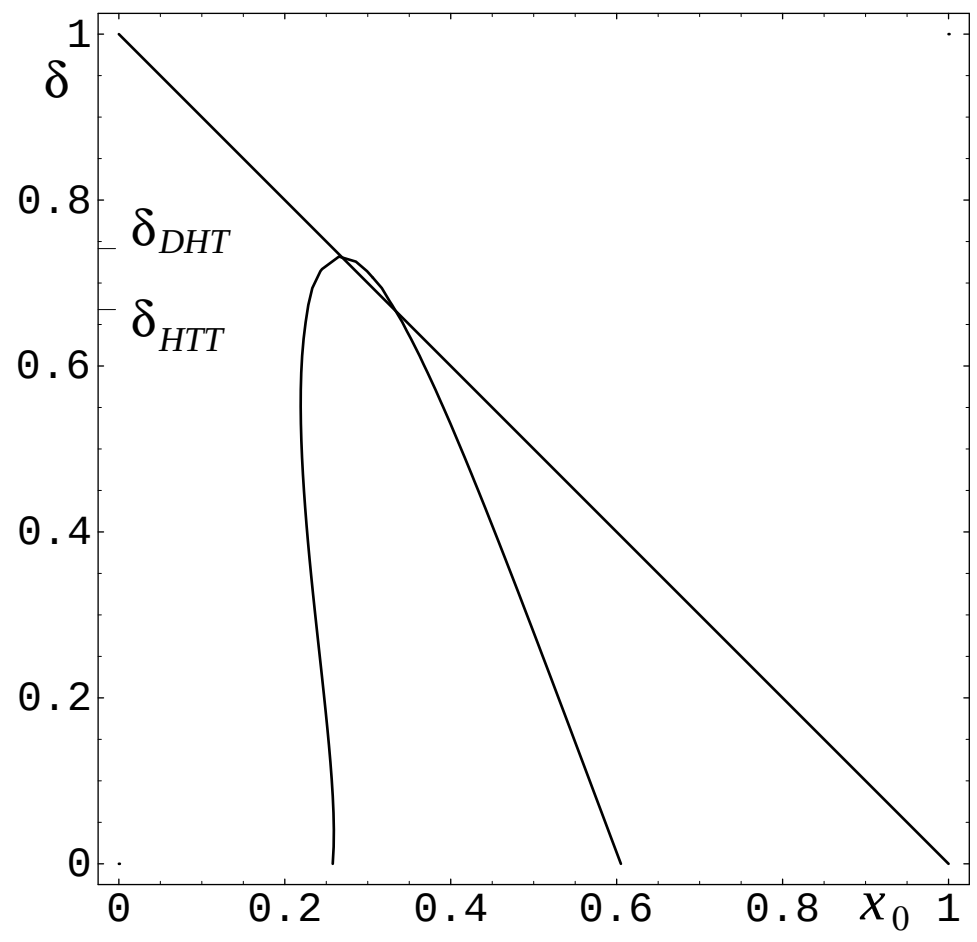

Fig. 5.

Bifurcations of codim 3 may be detected checking the non-degenerate conditions for bifurcation of codim 2.

In [2] it is proved that for parameters

$$a = 1/27, \ \beta = 9/4, \ \delta = 2/3$$

in the system there is a bifurcation "the focus singularity with nilpotent linear part" ( $HTT$).

It is proved also in this paper that for

$$a = -45 + 26\sqrt{3}, \ \beta = (5 + 3\sqrt{3})/4, \ \delta = \sqrt{3} - 1$$

in the sysyem there is a bifurcation "the cusp singularity with nilpotent linear part" ( $DTT$).

Parametric neighborhoods of listed local bifurcations are well Known and we know about sets of non-local bifurcations there. For example, with bifurcations $HT$ it takes place non-local separatrix loop bifurcation $L_{\mathrm{l}}\,(L_{\mathrm{r}})$, with bifurcation $HTT$ it takes place non-local double cycle bifurcation $C$ and so on.

So we have information about the existence and construction of the non-local bifurcation surfaces in regions of the parameter space near sets of local bifurcations. More detailed study of the surfaces may be fulfilled with the help of the qualitative theory of differential equations and numerical tests.

As a result of our investigation we obtain a global bifurcation diagram.

We do not describe here the diagram and we do not describe the process of its construction because this description repeats analogous investigations from other papers (see [7, 11, 12] and references therein).

To give an idea of the construction of our bifurcation diagram and to summarize our results we give the "most complex" section of the bifurcation diagram by the plane $\delta = c$ for value of $c$ near 0) (Fig. 6).

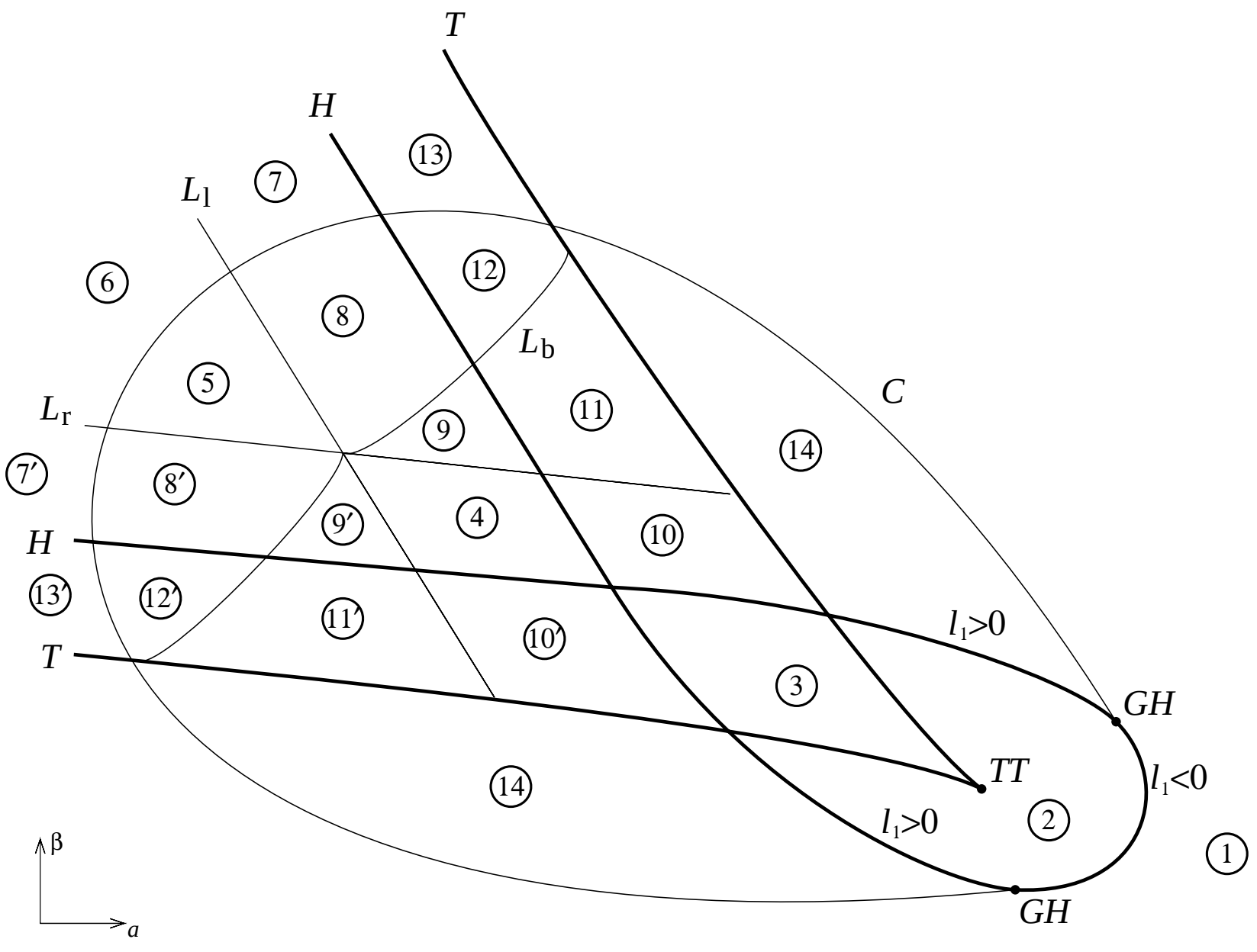

Fig. 6.

In the section the parameter space is $R^2_+$ with the coordinates $a, \beta$, codim one bifurcation sets are curves and codim two bifurcation sets are points with the same notations as corresponding bifurcations.

In our figure we sketch only disposition of bifurcation curves and points without following the scale because they are very close to each other sometimes and it is not possible to draw their full-size images.

In Fig. 7 are given schemes of rough phase portraits of system (2).

We show only trajectories within the set $\Pi$ and we do not sketch the boundary steady states. Stable steady states are depicted by black circles, ubstable steady states are depicted by white circles, saddles are depicted by half-black circles. Unstable limit cycles is depicted by dashed curves.

We do not show the phase portrait for regions $7'$–$14'$ of the bifurcation diagram. These portraits may be obtained from the portraits for regions 7–14 by interchanging two steady state which are non-saddles.

To complete the list of all phase portraits we show schemes of trajectories for two regions 15–16 of the global bifurcation diagram which are absent in the section.

We are thankful to professor Dongmei Xiao (Shanghai Jiao Tong University) for discussions and remarks.

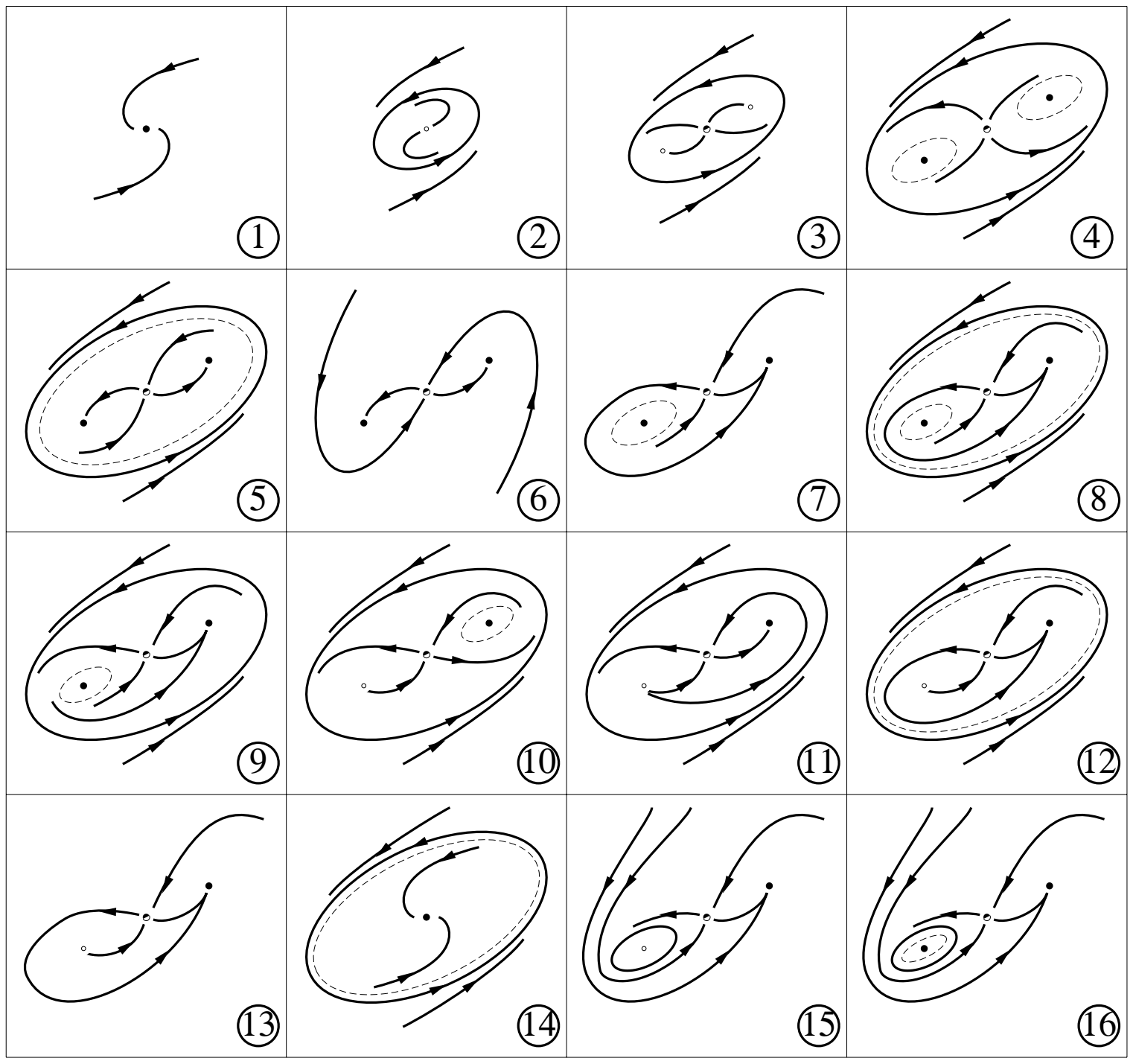

Fig. 7.

Sobolev Institute of Mathematics,
Novosibirsk
E-mail: treskov@math.nsc.ru volok@math.nsc.ru